\documentclass[12pt,a4paper]{amsart}
\usepackage{amsfonts}
\usepackage{amsthm}
\usepackage{amsmath}

\usepackage{amscd}
\usepackage{t1enc}
\usepackage[mathscr]{eucal}
\usepackage{indentfirst}
\usepackage{graphicx}
\usepackage{graphics}
\usepackage{pict2e}
\usepackage{enumerate}
\usepackage{array}
\usepackage{blindtext}
\usepackage[dvipsnames]{xcolor}
\definecolor{bred}{rgb}{0.8,0,0}
\usepackage{hyperref}
\usepackage{cleveref}
\hypersetup{colorlinks,linkcolor={blue},citecolor={bred},urlcolor={blue}}
\usepackage{epic}
\numberwithin{equation}{section}
\usepackage[margin=2.9cm]{geometry}
\usepackage{epstopdf}
\usepackage{mathtools}
\graphicspath{ {./PhD Images/} }

\Crefname{assumptionH}{\textbf{H}\hspace{-3pt}}{\textbf{H}\hspace{-3pt}}
\crefname{assumptionH}{\textbf{H}}{\textbf{H}}

\newtheorem{assumptionA}{\textbf{A}\hspace{-3pt}}
\Crefname{assumptionA}{\textbf{A}\hspace{-3pt}}{\textbf{H}\hspace{-3pt}}
\crefname{assumptionA}{\textbf{A}}{\textbf{A}}

\usepackage{amsmath,amssymb}

\theoremstyle{plain}
\newtheorem{Th}{Theorem}[section]
\newtheorem{Lemma}[Th]{Lemma}
\newtheorem{Cor}[Th]{Corollary}
\newtheorem{Prop}[Th]{Proposition}

\theoremstyle{definition}

\newtheorem{Rem}[Th]{Remark}
\newtheorem{?}[Th]{Problem}

\newtheorem{Note}[Th]{Note}

\author{Tim Johnston, Sotirios Sabanis}

\title{A Strongly Monotonic Polygonal Euler Scheme}

\begin{document}

	\begin{abstract}
In recent years tamed schemes have become an important technique for simulating SDEs and SPDEs whose continuous coefficients display superlinear growth. The taming method, which involves curbing the growth of the coefficients as a function of stepsize, has so far however not been adapted to preserve the monotonicity of the coefficients. This has arisen as an issue particularly in \cite{articletam}, where the lack of a strongly monotonic tamed scheme forces strong conditions on the setting.

In the present work we give a novel and explicit method for truncating monotonic functions in separable Hilbert spaces, and show how this can be used to define a polygonal (tamed) Euler scheme on finite dimensional space, preserving the monotonicity of the drift coefficient. This new method of truncation is well-defined with almost no assumptions and, unlike the well-known Moreau-Yosida regularisation, does not require an optimisation problem to be solved at each evaluation. Our construction is the first infinite dimensional method for truncating monotone functions that we are aware of, as well as the first explicit method in any number of dimensions.
	\end{abstract}

	\maketitle
\section{Introduction}
Whilst the Euler scheme method for approximating SDEs is known to be very robust in the sense of in probability convergence, see \cite{113194bd7a444bad9e31624261c18d4b}, it has also been known for almost fifteen years that in the case of superlinear coefficients such Euler schemes can possess extremely bad convergence properties in $L^p$, see \cite{Hutzenthaler_2010}. In fact the moments of the scheme itself are known to be unbounded for such equations as the stepsize tends to zero, as well as the $L^p$ difference between the scheme and the true solution.

This discovery led to the development of so-called `tamed' schemes firstly in \cite{Hutzenthaler_2012, Hutzenthaler_2015}, and subsequently in papers such as \cite{Sabanis_2013}, following Krylov's approach to Euler polygonal schemes, and using different methods in \cite{https://doi.org/10.48550/arxiv.1102.0662, Hutzenthaler_2017}. Also, related truncated Euler-Maruyama schemes appeared in \cite{MAO2015370, MAO2016362}. The idea of these schemes is that curbing the magnitude of the drift (and in later papers like \cite{kumar2016milstein, Sabanis_2016} the diffusion) coefficient of the SDE as a function of stepsize allows one to control the moments of the scheme (often via a one-sided Lipschitz condition), which then allows for control of the $L^p$ difference between the scheme and the true solution. This idea has also been fruitful for the development of stochastic algorithms such as \cite{tula, https://doi.org/10.48550/arxiv.2006.14514}.

Similar ideas have been introduced in \cite{articletam,Hutzenthaler_2020, Wang_2020} in the context of numerical methods for non-linear SPDEs. In the first of these articles the fact that the taming method did not preserve the monotonicity condition meant that strong conditions had to be assumed of the setting, in particular the decomposition of the drift into the sum of a superlinear and linearly growing component each defined on different spaces.

In this paper we construct an explicit truncation method for monotonic functions on finite and infinite dimensional separable Hilbert spaces, and demonstrate how this can be used to define a tamed Euler scheme on $\mathbb{R}^d$ that preserves the monotonicity of the drift coefficient. Additionally we show how a certain natural decomposition of so-called `strongly monotonic functions' can be used to extend our construction to functions of this type. In order to convey that our scheme preserves a richer class of properties than other tamed schemes, we choose to denote it as an example of a broader class of `polygonal Euler schemes', of which any Euler scheme for which the coefficients depend on the stepsize is an example.

Previously, a non-explicit construction of such a truncation method was given in Lemma 3 in \cite{doi:10.1137/1129033}, under stronger conditions and in finite dimensional space. To our knowledge this is the first infinite dimensional truncation procedure for monotonic functions preserving the monotonicity condition that has appeared in the literature, and the first explicit procedure in any number of dimensions. In addition our truncation procedure yields a Lipschitz function in the case where the original function is Lipschitz inside the radius of truncation.

\section{Truncating Monotonic Functions}
Our construction is given in the separable Hilbert space setting, which of course can be taken to be either finite or infinite dimensional. Let $\mathcal{H}$ be a separable Hilbert space with inner product $\langle \cdot, \cdot \rangle$ and norm $\lvert \cdot \rvert$. Let $B_r(x) \subset \mathcal{H}$ denote the open ball of radius $r$ around $x \in \mathcal{H}$. If $f: \mathcal{H} \to \mathcal{H}$ is a function let
\[
Lip_A(f):=\sup_{x,y \in A, x\not =y} \frac{\lvert f(x)-f(y)\rvert}{\lvert x-y\rvert},
\]
and let us call a function for which $Lip_A(f)$ is finite `Lipschitz on $A$'. By $Lip(f)$ we shall denote the global Lipschitz constant of $f$. A function $f$ satisfying
\begin{equation}\label{eq: mon cond}
\langle f(x)-f(y) , x-y \rangle \geq 0, \;\; x,y \in \mathcal{H},
\end{equation}
shall be called a `monotonic function'. Note that we have denoted mappings \linebreak $f:\mathcal{H}\to\mathcal{H}$ as `functions' so as to unify terminology in the finite and infinite dimensional case. In Theorem \ref{th: mon} below we present our construction for a truncation procedure preserving the monotonicity condition.

\begin{Th} \label{th: mon}
Let $f: \mathcal{H} \to \mathcal{H}$ satisfy the monotonicity property \eqref{eq: mon cond} and suppose
\begin{equation}\label{eq: f bound}
\sup_{x \in B_r(x_0)} \lvert f(x) \rvert \leq R,
\end{equation}
holds for some $x_0 \in \mathcal{H}$, $r\geq 2, R>0$. Then there exists a truncated function $\tilde{f}_{R, r, x_0}: \mathcal{H} \to \mathcal{H}$ such that
	\begin{enumerate}[(i)]
			\item For every $x \in B_{r-2}(x_0)$, $\tilde{f}_{R, r, x_0}(x)=f(x)$
		\item The function $\tilde{f}_{R, r, x_0}$ obeys the bound
	\begin{equation}
	\lvert \tilde{f}_{R, r, x_0}(x) \rvert \leq R (1+\min\{r, \lvert x-x_0 \rvert \}). \nonumber
	\end{equation}
			\item If $f$ is Lipschitz on $B_r$, then $\tilde{f}_{R, r, x_0}$ is globally Lipschitz and in particular
	\begin{equation}
	Lip(\tilde{f}_{R, r, x_0}) \leq Lip_{B_r}(f) +Rr. \nonumber
	\end{equation}
		\item The monotonicity property of $f$ is conserved, i.e. $\tilde{f}$ satisfies
		\begin{equation}
	\langle \tilde{f}_{R, r, x_0}(x)-\tilde{f}_{R, r, x_0}(y) , x-y \rangle  \geq 0, \;\; x,y \in \mathcal{H}. \nonumber
		\end{equation}
	\end{enumerate}
Specifically, $f$ is given as
\[
\tilde{f}_{R, r, x_0}(x) := t(x)f(x) + R s ( x )(x-x_0) ,
\]
where\textbf{}
\begin{equation*}
t(x) := \begin{cases}
1, & \lvert x -x_0 \rvert \leq r-1 \\
r-\lvert x -x_0 \rvert, & r-1 < \lvert x -x_0 \rvert < r \\
0, & \lvert x -x_0 \rvert \geq r
\end{cases} \nonumber , \;\;\;
\end{equation*}
\begin{equation*}
s(x) := \begin{cases}
0, & \lvert x -x_0 \rvert \leq r-2\\
\lvert x -x_0 \rvert -r+2, & r-2 < \lvert x -x_0 \rvert < r-1 \\
1, &   r-1 \leq \lvert x -x_0  \rvert \leq r \\
\frac{ r }{\lvert x -x_0  \rvert }, &  \lvert x -x_0 \rvert \geq r
\end{cases} . \nonumber
\end{equation*}
\end{Th}

\begin{figure}[h] \label{fig 1}
	\caption{The construction of $\tilde{f}_{R, r, x_0}$}
	\centering
	\includegraphics[width=14cm, height=8.9cm]{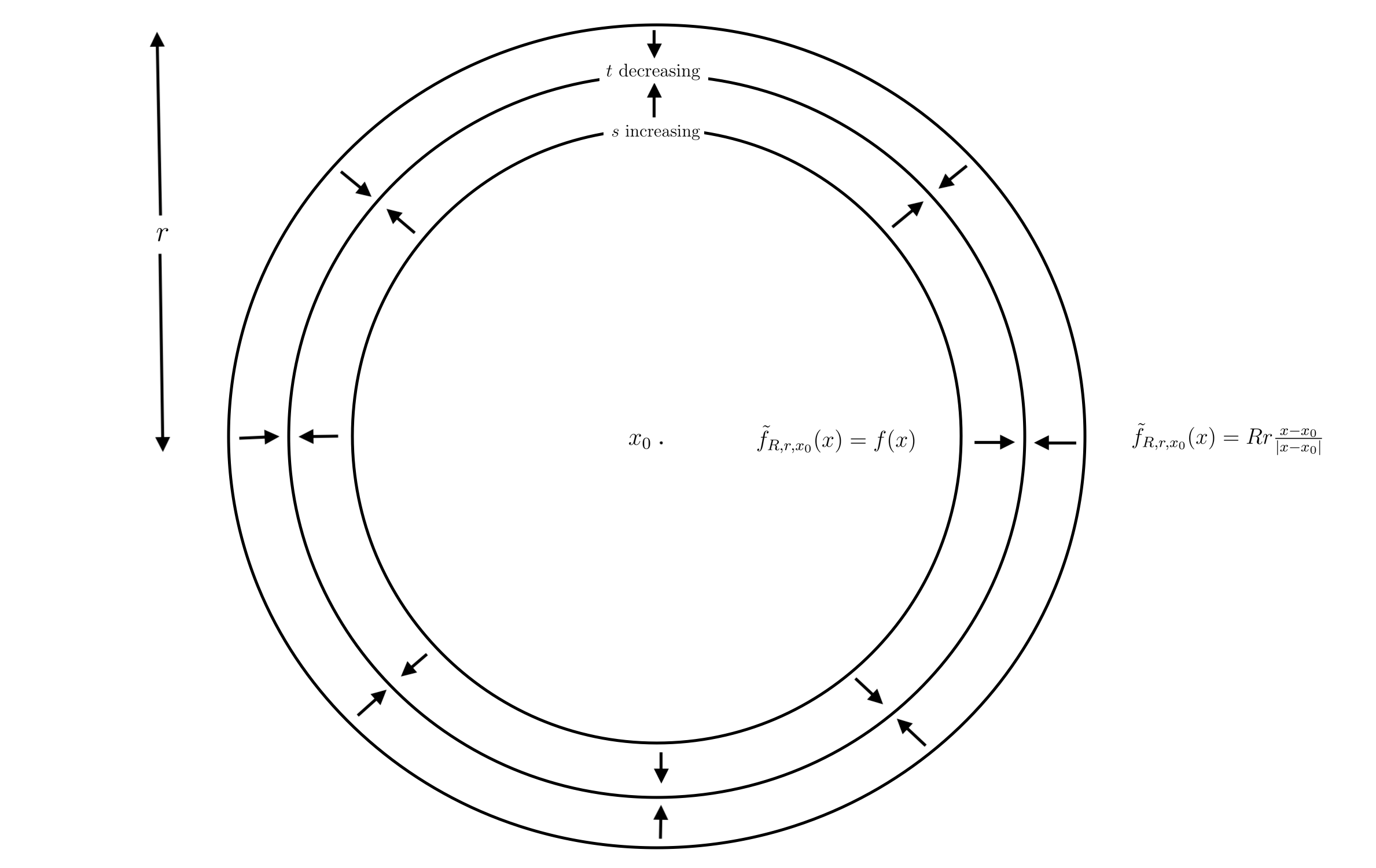}.
\end{figure}
\begin{proof}

To ease notation we write $\tilde{f}$ for $\tilde{f}_{R, r, x_0}$. From the definition of $\tilde{f}$ it follows that $t=1$ and $s=0$ on $B_{r-2}(x_0)$ and therefore we obtain (i) immediately. Furthermore (ii) follows from \eqref{eq: f bound} and the triangle inequality since $\lvert t(x) f(x) \rvert \leq R$ and \linebreak $\lvert Rs(x)(x-x_0) \rvert \leq R\min\{r, \lvert x-x_0 \rvert\}$.

To prove (iii) first recall that for $A \subset \mathcal{H}$ and Lipschitz functions $f:A \to \mathcal{H}$, $g: A \to \mathbb{R}$ one has the bound
\begin{equation}\label{eq: Lip bd}
Lip_A(fg)\leq \sup_{x\in A}\lvert f(x) \rvert Lip_A(g)+\sup_{x\in A}\lvert g(x) \rvert Lip_A(f).
\end{equation}
Now let us define
\begin{equation*}
\Omega_1:= \bar{B}_{r-2}, \;\; \Omega_2 := \bar{B}_{r-1} \setminus B_{r-2}, \;\; \Omega_3 :=  \bar{B}_{r} \setminus B_{r-1}, \;\; \Omega_4 := \mathcal{H} \setminus B_{r},
\end{equation*}
corresponding to (the completion of) the four region in Figure 1. Then if $x,y \in \Omega_1$ one has
\begin{equation*}
\lvert \tilde{f}(x)-\tilde{f}(y) \rvert \leq Lip_{B_r}(f)\lvert x-y\rvert,
\end{equation*}
since $f=\tilde{f}$ on $B_{r-2}$. For $x,y \in \Omega_2$, noting that $Lip(s)=1$ and $t=1$ on $\Omega_2$ one uses \eqref{eq: Lip bd} to obtain
\begin{equation*}
\lvert \tilde{f}(x)-\tilde{f}(y) \rvert \leq (Lip_{B_r}(f)+Rr)\lvert x-y\rvert.
\end{equation*}
Similarly, for $x,y \in \Omega_3$ since $Lip(t)=1$ and $\sup_{x \in \Omega_3}\lvert f(x) \rvert \leq R$  by \eqref{eq: f bound}, one uses \eqref{eq: Lip bd} again to obtain
\begin{equation*}
\lvert \tilde{f}(x)-\tilde{f}(y) \rvert  \leq (Lip_{B_r}(f) +2R)\lvert x-y \rvert.
\end{equation*}
Finally for $x,y \in \Omega_4$, since $\lvert x-x_0\rvert, \lvert y-x_0\rvert \geq r$, and therefore $\frac{r^2}{\lvert x-x_0\rvert\lvert y-x_0\rvert }-1<0$, one calculates
\begin{align}
R^2\lvert x-y\rvert^2&-\lvert \tilde{f}(x)-\tilde{f}(y)\rvert^2 =R^2\lvert (x-x_0)-(y-x_0)\rvert^2-\biggr \lvert \frac{ Rr }{\lvert x -x_0  \rvert }(x-x_0)-\frac{ Rr }{\lvert y -x_0  \rvert }(y-x_0) \biggr \rvert^2 \nonumber \\
& =R^2\biggr(\lvert x-x_0\rvert^2+\lvert y-x_0\rvert^2-2r^2+2\langle x-x_0, y-x_0 \rangle \biggr [\frac{r^2}{\lvert x-x_0\rvert\lvert y-x_0\rvert }-1\biggr] \biggr) \nonumber \\
&\geq R^2(\lvert x-x_0\rvert^2+\lvert y-x_0\rvert^2-2r^2+2r^2-2\lvert x-x_0\rvert\lvert y-x_0\rvert)
\nonumber \\
& \geq R^2(\lvert x-x_0\rvert-\lvert y-x_0\rvert)^2 \geq 0,
\end{align}
so that $\lvert  \tilde{f}(x)-\tilde{f}(y)\rvert \leq R\lvert x-y\rvert$. Therefore one obtains (iii) locally on each $\Omega_i$, that is, since $r \geq 2$
\begin{equation}\label{eq: Lipschitz bd}
\lvert  \tilde{f}(x)-\tilde{f}(y)\rvert \leq (Lip_{B_r}(f)+Rr) \lvert x-y \rvert, \;\;x,y \in \Omega_i,
\end{equation}
holds for $i=1,2,3,4$. Now if we let $x,y \in \mathcal{H}$ be arbitrary, we see that \eqref{eq: Lipschitz bd} continues to hold. Indeed, let $l \subset \mathcal{H}$ denote the line segment connecting $x$ and $y$. Since every straight line intersects the boundary of any ball at most twice, there exist at most $6$ points where $l$ intersects any $\partial \Omega_i$. Let us denote these points in order from $x$ to $y$ as $l_1, ... , l_{k-1} \in l$ with $l_0=x, l_k=y$. Then since each $\Omega_i$ is a closed subset of $\mathcal{H}$ and $\cup_{i=1}^4 \Omega_i = \mathcal{H}$, the line segment connecting $l_i$ and $l_{i+1}$ for any $i=0, ..., k-1$ is entirely contained in some $\Omega_i$, otherwise there would have to be another point of intersection with some $\partial \Omega_i$ between them. Therefore since the $\Omega_i$ are closed
\begin{align}
\lvert  \tilde{f}(x)-\tilde{f}(y)\rvert &\leq \sum_{i=1}^k \lvert  \tilde{f}(l_i)-\tilde{f}(l_{i-1})\rvert \nonumber \\
& \leq \sum_{i=1}^k (Lip_{B_r}(f)+Rr) \lvert l_i-l_{i-1} \rvert\nonumber \\
& = (Lip_{B_r}(f)+Rr) \lvert x-y \rvert.
\end{align}
It remains to prove that $\tilde{f}$ preserves the monotonicity property. Let us first show that $\tilde{f}$ obeys the monotonicity condition piecewise on each $\Omega_i$, i.e.
		\begin{equation*}
	\langle \tilde{f}(x)-\tilde{f}(y), x-y \rangle \geq 0, \;\; x,y \in \Omega_i ,
	\end{equation*}
	for $i=1,2,3,4$. By construction $\tilde{f}=f$ on $\Omega_1$, so clearly $\tilde{f}$ obeys the monotonicity condition on $\Omega_1$. Now assume without loss of generality that $x,y \in \Omega_2$ satisfies $\lvert x-x_0 \rvert  \geq \lvert y-x_0 \rvert$, so that $s(x) \geq s(y) $ and therefore
	\begin{align}
	\langle \tilde{f}(x) - &\tilde{f}(y), x-y \rangle
	= \langle f( x) - f(y), x-y \rangle
\nonumber \\
&	+ R\langle s( x ) (x-x_0) - s( y ) (y-x_0), (x-x_0)-(y-x_0)\rangle \nonumber \\
& \geq s(x) \lvert x-x_0\rvert^2+s(y) \lvert y-x_0\rvert^2 - (s(x)+s(y)) \langle x-x_0, (x-x_0)-(y-x_0)\rangle \nonumber \\
	& \geq R(s( x ) \lvert x-x_0 \rvert -s(y)\lvert y-x_0 \rvert)(\lvert x-x_0 \rvert-\lvert y-x_0 \rvert) \geq 0.
	\end{align}
	For $x, y \in \Omega_3$ also satisfying $\lvert x-x_0 \rvert  \geq \lvert y-x_0 \rvert$, since $ \lvert \lvert x -x_0\rvert - \lvert y-x_0 \rvert \rvert \leq   \lvert x-y \rvert$, and by \eqref{eq: f bound} and the monotonicity of $f$, one has
	\begin{align*}
	\langle \tilde{f}(x) - \tilde{f}(y), x-y \rangle&= R\lvert x-y \rvert^2 + t(y) \langle f(x)-f(y), x-y \rangle
	+  (t(x)- t(y) ) \langle f(x), x-y \rangle   \\
	& \geq R \lvert x-y \rvert^2 - \lvert \lvert x -x_0\rvert - \lvert y-x_0 \rvert \rvert\lvert f(x) \rvert \lvert x-y \rvert  \\
	& \geq (R - \lvert f(x) \rvert) \lvert x-y \rvert^2 \geq 0 .
	\end{align*}
	Finally for $x, y \in \Omega_4$, using Cauchy-Schwarz
	\begin{align*}
	\langle \tilde{f}(x) - \tilde{f}(y), x-y \rangle&= R  r \biggr  \langle \frac{x-x_0}{\lvert x -x_0\rvert}-\frac{y-x_0}{\lvert y-x_0 \rvert}, (x-x_0)-(y-x_0) \biggr \rangle \geq 0.
	\end{align*}
Now let us show that this implies that $\tilde{f}$ is monotonic on the whole of $\mathcal{H}$. Let $x, y \in \mathcal{H}$ be arbitrary. As before let $l \subset \mathcal{H}$ denote the line segment connecting $x$ and $y$ and $l_1,..., l_{k-1}$ the intersections with $\cup_{i=1}^4 \partial \Omega_i$ from $x$ to $y$. Thus $\langle f( l_{i})-f(l_{i-1}), l_{i}-l_{i-1} \rangle \geq 0$ for $i=0, ..., k$, and therefore since $l_1-l_0, l_2-l_1, ..., l_{k}-l_{k-1}$ are each collinear with $x-y$, one has since the $\Omega_i$ are closed that
	\begin{align}
\langle \tilde{f}(x)-\tilde{f}(y), x-y \rangle = \sum_{i=1}^{k} 	\langle f( l_{i})-f(l_{i-1}), l_{i}-l_{i-1} \rangle \frac{\lvert x-y \rvert}{\lvert   l_{i}-l_{i-1}\rvert} \geq 0, \nonumber
\end{align}
as required.
\end{proof}
We can now adapt Theorem \ref{th: mon}  to the setting where $f$ obeys the strong monotonicity assumption
\begin{equation} \label{eq: strong mon}
\langle f(x)-f(y) , x-y \rangle \geq L\lvert x-y\rvert^2 \;\; x,y \in \mathcal{H} .
\end{equation}
for some $L>0$. A function obeying this property can always be decomposed as the sum of a linear function and a monotonic function, that is $f:\mathcal{H}\to\mathcal{H}$ satisfying \eqref{eq: strong mon} can always be written as $f(x)=g(x)+Lx$ for a function $g:\mathcal{H}\to\mathcal{H}$ obeying the monotonicity condition \eqref{eq: mon cond}. Since such a function $f$ can never be bounded, we shall instead apply Theorem \ref{th: mon} to the component $g$, in order to produce a function that is the sum of a bounded and linear part.
\begin{Cor} \label{cor: strong mon}
 Let $f: \mathcal{H} \to \mathcal{H}$ obey the strong monotonicity condition \eqref{eq: strong mon} with constant $L>0$. Let $g:\mathcal{H}\to\mathcal{H}$ be such that $f(x)=g(x)+Lx$ (so $g(x):=f(x)-Lx$). Suppose that
	\begin{equation}\label{eq: g bound}
	\sup_{x \in B_r(x_0)} \lvert g(x)\rvert \leq R,
	\end{equation}
for some $x_0 \in \mathcal{H}$, $r\geq 2, R>0$. Then there exists a function $\bar{f}_{R, r, x_0}:\mathcal{H}\to\mathcal{H}$ such that
	\begin{enumerate}[(i)]
		\item For every $x \in B_{r-2}$, $\bar{f}_{R, r, x_0}(x)=f(x)$
		\item The function $\bar{f}_{R, r, x_0}$ obeys the bound
		\begin{equation*}
		\lvert \bar{f}_{R, r, x_0}(x) \rvert \leq R (1+\min\{r, \lvert x-x_0 \rvert \})+L\lvert x \rvert .
		\end{equation*}
					\item If $f$ is Lipschitz on $B_r$, then $\bar{f}_{R, r, x_0}$ is globally Lipschitz and one has the bound
		\begin{equation*}
		Lip(\bar{f}_{R, r, x_0}) \leq Lip_{B_r}(g) +Rr+L .
		\end{equation*}
		\item The strong monotonicity property of $f$ is conserved, i.e. $\bar{f}_{R, r, x_0}$ satisfies
		\begin{equation*}
		\langle \bar{f}_{R, r, x_0}(x)-\bar{f}_{R, r, x_0}(y) , x-y \rangle  \geq L\lvert x-y \rvert^2, \;\; x,y \in \mathcal{H} .
		\end{equation*}
	\end{enumerate}
	Specifically, $f$ is given as
	\begin{equation*}
\bar{f}_{R, r, x_0}(x) := \tilde{g}_{R, r, x_0}(x)+Lx ,
	\end{equation*}
	for $\tilde{g}_{R, r, x_0}$ as in Theorem \ref{th: mon}.
\end{Cor}
\begin{proof}
Observe that by the strong monotonicity property \eqref{eq: strong mon} for $f$, the function $g$ is monotonic. Therefore one can apply Theorem \ref{th: mon} to $g$ to produce the function $\tilde{g}_{R, r, x_0}$ that is also monotonic. Then (i), (ii) and (iii) follow immediately from Theorem \ref{th: mon} and the definition of $\bar{f}_{R, r, x_0}$. Furthermore, we see that for every $x, y \in \mathcal{H}$
\begin{equation*}
	\langle \bar{f}_{R, r, x_0}(x)-\bar{f}_{R, r, x_0}(y) , x-y \rangle =  	\langle \tilde{g}_{R, r, x_0}(x)-\tilde{g}_{R, r, x_0}(y) , x-y \rangle +L\lvert x-y \rvert^2 \geq L\lvert x-y \rvert^2,
\end{equation*}
proving (iv).
\end{proof}

\begin{Rem} \label{rem: comparison}
In comparison to the well-known Moreau-Yosida regularisation (which requires an optimisation problem to be solved for each evaluation), our approach is explicit and therefore straightforward to evaluate. However unlike the Moreau-Yosida regularisation it does not preserve minimisers unless the minimal point lies inside the ball of truncation, and it does not yield a Lipschitz function unless the original function is Lipschitz inside the ball of truncation.
\end{Rem}

\section{A Strongly Monotonic Euler Scheme}
In this section we adopt the use Theorem \ref{th: mon} in the setting of polygonal Euler approximations following Krylov's paradigm, see for example \cite{Krylov:90}, where the coefficients of the approximate scheme dependent directly on the step size. In our case the sequence of coefficients is given in Proposition \ref{th: strong mon} below via Theorem \ref{th: mon}, in such a way as to preserve the strong monotonicity condition, whilst also being sufficiently bounded so as to prove $L^p$ convergence of the scheme in the presence of superlinear coefficients.

Specifically we let $b:\mathbb{R}^d \to \mathbb{R}^d$ and $\sigma:\mathbb{R}^d \to \mathbb{R}^{d \times d}$ be functions, and $\eta$ an $\mathbb{R}^d$-valued random variable independent of a $d$-dimensional Wiener martingale $(W_t)_{t\geq 0}$. Consider the stochastic differential equation
\begin{equation} \label{eq: mainSDE}
dX(t) = b(X(t)) dt+\sigma(X(t)) dW_t, \;\;\; X(0) = \eta,  \;\; t \in [0,T].
\end{equation}
and consider furthermore an Euler-Krylov approximation to \eqref{eq: mainSDE} of the form
\begin{equation} \label{eq: monotonic-polygonal-scheme}
dX_n(t) = b_n(X_n(\kappa_n(t)))dt+\sigma(X_n(\kappa_n(t))) dW_t, \;\;\; X_n(0) = \eta, \;\;\; t \in[0,T] ,
\end{equation}
where $\kappa_n(t) :=  [nt]/n$ is the projection backwards to the grid of width $1/n$ and $b_n: \mathbb{R}^d \to \mathbb{R}^d$ are the sequence of functions given below in Proposition \ref{th: strong mon}. Denoting the inner product of $x,y \in \mathbb{R}^d$ by $xy$, we  assume of $b, \sigma$ that

\begin{assumptionA} \label{assump: mon, lip} There exists a constant $c>0$ such that
\begin{equation} \label{eq: strong-monotonicity}
(b(x) - b(y))(x-y) \leq -L \lvert x - y \rvert^2, \qquad\forall x,y \in \mathbb{R}^d, \nonumber
\end{equation}
\begin{equation}
\lvert \sigma(x)-\sigma(y) \rvert \leq c \lvert x-y \rvert ,\qquad\forall  x \in \mathbb{R}^d. \nonumber
\end{equation}
\end{assumptionA}

\begin{assumptionA}
\label{assump: in cond} There exists a positive constant $p_0>0$ such that $E \lvert \eta \rvert^{p_0} < \infty $.
\end{assumptionA}

\begin{assumptionA}
\label{assump: poly growth} There exist positive constants $H,l>0$ such that
\begin{equation*}
\lvert g(x)-g(y) \rvert \leq H(1+\lvert x\rvert+\lvert y \rvert)^l \lvert x-y \rvert, \;\;\; x,y \in \mathbb{R}^d.
\end{equation*}
where $g$ is given in \eqref{eq: b decomp}.
\end{assumptionA}

Note that A\ref{assump: mon, lip} implies that the function $-b$ obeys the strong monotonicity property with constant $L>0$. In particular we can write
\begin{equation}\label{eq: b decomp}
b(x)=b(0)-g(x)-Lx,
\end{equation}
for a monotonic function $g: \mathbb{R}^d\to\mathbb{R}^d$, as in Corollary \ref{cor: strong mon} (writing $b(0)-g(x)$ for the `monotonic' part of $b$ will be marginally more convenient
 in the following arguments since this implies $g(0)=0$).
 \begin{Note}
 From now on we let $C>0$ (and occasionally $C_1, C_2,... >0$) denote generic positive constants independent of $n$ and $t$ that changes from line to line.
 \end{Note}

\begin{Prop} \label{th: strong mon}
	Let A\ref{assump: mon, lip}, A\ref{assump: in cond} and A\ref{assump: poly growth} hold. Then there exists $m\in \mathbb{N}$ and for $n\geq m$ functions $b_n:\mathbb{R}^d\to \mathbb{R}^d$ (given explicitly in \eqref{eq: b_n def} below) such that
	\begin{enumerate}[(i)]
		\item The strong monotonicity property of $b$ is conserved in $b_n$, i.e.
		\begin{equation*}
		(b_n(x)-b_n(y))(x-y) \leq -L\lvert x-y\rvert^2, \;\; x,y \in \mathbb{R}^d ,
		\end{equation*}
		\item The functions $b_n$ each obey the bound
		\begin{equation*}
		\lvert b_n(x) \rvert \leq L\lvert x \rvert+n^{1/2} (1+\min \{\lvert x \rvert ,s_{n^{1/2}} \}) ,
		\end{equation*}
		for $s_n = H^{-\frac{1}{l+1}}n^{\frac{1 }{2(l+1)}}-1$.
		\item For every $\mathbb{R}^d$-valued random variable $X$ satisfying $E \lvert X\rvert^{p_0}  <\infty$ there exist a constant $C(E \lvert X\rvert^{p_0})>0$ such that
		\begin{equation*}
		E \lvert b(X)-b_n(X) \rvert^p \leq Cn^{\frac{(l+2)p-p_0}{2(l+1)} } E \lvert X\rvert^{p_0}.
		\end{equation*}
	\end{enumerate}
\end{Prop}

\begin{proof}
Let $g$ be as in \eqref{eq: b decomp}. Observe by the smoothness assumption A\ref{assump: poly growth} there exists $m \in \mathbb{N}$ and a sequence $(s_n)_{n \geq m} \in \mathbb{R}^d$ such that $s_n \to \infty$ and
\begin{equation*}
\sup_{x \in B_{s_n}(0)} \lvert g(x) \rvert \leq n^{1/2}, \;\; n \geq m.
\end{equation*}
Furthermore, since $g(0)=0$ one has
\begin{equation*}
\lvert g(x) \rvert \leq H(1+\lvert x \rvert)^l\lvert x \rvert \leq H(1+\lvert x \rvert)^{l+1},
\end{equation*}
so one can choose $s_n = H^{-\frac{1}{l+1}}n^{\frac{1 }{2(l+1)}}-1$ and $m$ sufficiently large that $s_n \geq 2$ for $n \geq m$. Following from \eqref{eq: b decomp}, similarly to as in Corollary \ref{cor: strong mon}, we define
\begin{equation}\label{eq: b_n def}
b_n(x):=b(0)-Lx-\widetilde{g}_{n^{1/2}, s_n,0}(x)
\end{equation}
for $\widetilde{g}_{n^{1/2}, s_n,0}$ as given in Theorem \ref{th: mon}. From Theorem \ref{th: mon} and following the argument of Corollary \ref{cor: strong mon} (with signs reversed) one immediately obtains (i) and (ii). To prove (iii) first recall that for a random variable $X$ on $\mathbb{R}^d$ such that $E\lvert X\rvert^{p_0}<\infty$, one has for every $q \in (0,p_0)$ that
	\begin{equation} \label{eq: key rv bd}
E \lvert X\rvert ^q1_{\{X>x \} } \leq x^{q-q_0}  E\lvert X\rvert^{q_0}  .
\end{equation}	
This inequality follows from Holder's and  Markov's inequality. Therefore, since $b$ and $b_n$ only differ outside the ball of radius $s_n-2$, one has by Theorem \ref{th: mon}, \eqref{eq: b decomp} and \eqref{eq: b_n def} that for an $\mathbb{R}^d$-valued random variable $X$
\begin{align*}
\lvert b(X)-b_n(X) \rvert ^p \leq C	(n^{ p/2} (1+\lvert X \rvert )^p+ \lvert g(X) \rvert^p) 1_{\{ \lvert X \rvert > s_n -2 \}}.
\end{align*}
Now note that
\begin{equation}\label{eq: s_n bound}
s_n \geq Cn^{\frac{1 }{2(l+1)}}.
\end{equation}
Indeed, \eqref{eq: s_n bound} clearly holds asymptotically, that is, there exists $N>0$ such that $s_n \geq Cn^{\frac{1 }{2(l+1)}}$ for $n\geq N$. One can then choose the constant $C>0$ sufficiently small that $s_n \geq Cn^{\frac{1 }{2(l+1)}}$ for all $m \leq n \leq N$.
Then, by \eqref{eq: key rv bd}, for all $n\geq m$ one has that
\begin{align}\label{Term 1}
n^{ p/2} E		 (1+\lvert X \rvert)^p 1_{\{\lvert X \rvert > s_n-2 \}} & \leq n^{ p/2} E(1+\lvert X \rvert)^p 1_{\{1+\lvert X \rvert > Cn^{\frac{1}{2(l+1)}} \}} \nonumber \\
& \leq Cn^{\frac{(l+2)p-p_0}{2(l+1)}} E	 (1+ \lvert X \rvert )^{p_0},
\end{align}
and by A\ref{assump: poly growth}
\begin{align}\label{Term 2}
E 	 \lvert g(X) \rvert^p 1_{\{\lvert X \rvert > s_n-2 \}}  &\leq C_1 E	(1+\lvert X\rvert^{p(l+1)}) 1_{\{	\lvert X\rvert > C_2n^{\frac{1}{2(l+1)}} \}} \nonumber \\
&\leq C_1 E(1+\lvert X \rvert)^{p(l+1)} 1_{\{	1+	\lvert X \rvert > C_2n^{\frac{1}{2(l+1)}} \}} \nonumber \\
& \leq Cn^{\frac{p(l+1)-p_0}{2(l+1)}} E(1+ \lvert X \rvert )^{p_0} .
\end{align}
Thus, if we assume $E\lvert X \rvert^{p_0} <\infty$, since $(1+x)^{p_0} \leq C(1+x^{p_0})$ and since the exponent in \eqref{Term 1} is more positive than the exponent in \eqref{Term 2}, the result follows.
\end{proof}

We now show that the scheme \eqref{eq: monotonic-polygonal-scheme} with $b_n$ as given in Proposition \ref{th: strong mon} converges to the true solution \eqref{eq: mainSDE}. For completeness we provide an entire proof.

\begin{Th} \label{thm:rate_of_convergence}
	Let A\ref{assump: mon, lip}, A\ref{assump: in cond} and A\ref{assump: poly growth} hold. Suppose $4(l+1) \leq p_0$, and consider $X_n$ as given in \eqref{eq: monotonic-polygonal-scheme}, for $b_n$ as defined in Proposition \ref{th: strong mon}. Then for every $p < p_0/(l+2)$ and every $n \geq m$
	\begin{equation}\label{eq: conv}
	E[\sup_{0\leq t \leq T}\lvert X(t) - X_n(t) \rvert^p] \leq Cn^{-r},
	\end{equation}
	where $r= \min \{\frac{1}{2(l+1)}(p_0(1 \vee 4/p)-(l+2)p), p/2 \}$.
\end{Th}

\begin{Rem}
One achieves an optimal $L^p$ convergence rate of $O(n^{-p/2})$ in Theorem \eqref{eq: conv} when $p_0$ in assumption A\ref{assump: in cond} is sufficiently large.
\end{Rem}

In order to prove Theorem \ref{thm:rate_of_convergence} we shall show that  $E\sup_{0 \leq t \leq T} \lvert X_n(t) \rvert^p$ is uniformly bounded in $n$. Key to these proofs shall be the property
\begin{align} \label{coercivity}
b(x)x \leq b(0)x-L\lvert x \rvert ^2\leq -\frac{L}{2} \lvert x \rvert ^2+ \frac{1}{2L} \lvert b(0) \rvert^2,
\end{align}
which follows from A\ref{assump: mon, lip}.
\begin{Prop}\label{thm: moment bounds}
	Let A\ref{assump: mon, lip}, A\ref{assump: in cond} and A\ref{assump: poly growth} hold. Suppose $4(l+1) \leq p_0$, and consider $X_n$ as given in \eqref{eq: monotonic-polygonal-scheme}, for $b_n$ as defined in Proposition \ref{th: strong mon}. Then, for every $0<p \leq p_0$ and every $n \geq m$
	\begin{equation*}
	\sup_{n\geq 1}	E[\sup_{0 \leq t \leq T} \lvert X_n(t) \rvert ^p] \leq C.
	\end{equation*}
\end{Prop}
\begin{proof}
	We use the notation $E^x$ to denote expectation with respect to the initial condition $X_n(0)=x \in \mathbb{R}^d$, or equivalently $E^x[\; \cdot\; ] = E[\; \cdot \; \vert X_n(0)=x]$. This is well defined since $X_n$ is the unique strong solution to \eqref{eq: monotonic-polygonal-scheme}. Let us first assume $p \geq 4$. By It{\^o}'s lemma and writing $b_n(X_n(\kappa_n(s)))X_n(s) = b_n(X_n(\kappa_n(s)))X_n(\kappa_n(s))+b_n(X_n(\kappa_n(s)))(X_n(s)-X_n(\kappa_n(s)))$ so as to apply \eqref{coercivity}, one obtains
	\begin{align}
	\lvert X_n(t) \rvert^2 \leq C \biggr( \lvert \eta \rvert ^2 &+\int^t_0[ b_n(X_n(\kappa_n(s)))X_n(s)+ \lvert \sigma(X_n(\kappa_n(s))) \rvert^2 ]ds + \int^t_0 \sigma(X_n(\kappa_n(s))) X_n(s) dW_s \biggr) \nonumber \\
	& \leq C \biggr (1+ \lvert \eta \rvert ^2 +\int^t_0 b_n(X_n(\kappa_n(s)))(X_n(s)-X_n(\kappa_n(s))) ds \nonumber \\
	& + \int^t_0 \lvert \sigma(X_n(\kappa_n(s))) \rvert^2 ds + \int^t_0 \sigma(X_n(\kappa_n(s))) X_n(s) dW_s \biggr) , \nonumber
	\end{align}
	at which point raising to the power of $p/2$ and using Holder's inequality yields
	\begin{align}
	\lvert X_n(t) \rvert^p \leq C \biggr (1+ \lvert \eta \rvert^p +&\int^t_0 \biggr( \lvert b_n(X_n(\kappa_n(s)))(X_n(s)-X_n(\kappa_n(s))) \rvert^{p/2}+ \lvert X_n(\kappa_n(s))\rvert^p \biggr) ds \nonumber \\
	& + \left( \int^t_0 \sigma(X_n(\kappa_n(s))) X_n(s) dW_s \right)^{p/2} \biggr), \nonumber
	\end{align}
	so that
	\begin{align} \label{eq:sup_estimate}
	\sup_{0 \leq t \leq T} \lvert X_n(t) \rvert^p \leq C \biggr (1+ &\lvert \eta \rvert^p +\int^T_0 \lvert \left ( b_n(X_n(\kappa_n(s)))(X_n(s)-X_n(\kappa_n(s))) \rvert^{p/2}+ \lvert X_n(\kappa_n(s))\rvert^p \right ) ds \nonumber \\
	& + \sup_{0 \leq t \leq T} \left( \int^t_0 \sigma(X_n(\kappa_n(s))) X_n(s) dW_s \right)^{p/2} \biggr).
	\end{align}
	Now one can observe that since by A\ref{assump: mon, lip} and Proposition \ref{th: strong mon} (ii), $b_n$ and $\sigma$ are both bounded by affine functions of $\lvert x \rvert$, it is a standard result that $E^x \lvert X_n(t) \rvert ^q < \infty$ for every $q>0$ and $n \in \mathbb{N}$ (see \cite{c711ce3d65fa4eba8f1294ddf46ec1d4}). Therefore the stochastic integral in the last term of \eqref{eq:sup_estimate} is a true martingale and one can apply the Burkholder-Davis-Gundy inequality to conclude
	\begin{align} \label{eq:expect_sup_estimate}
	E ^x\sup_{0 \leq t \leq T} \lvert X_n(t) \rvert^p &\leq C \biggr (1+ \lvert x \rvert^p +\int^T_0 E ^x \lvert b_n(X_n(\kappa_n(s)))(X_n(s)-X_n(\kappa_n(s))) \rvert^{p/2} ds \nonumber \\
	& +\int^T_0 E ^x\lvert X_n(\kappa_n(s))\rvert^p ds + E ^x \left( \int^T_0 \lvert \sigma(X_n(\kappa_n(s))) X_n(s)  \rvert^2 ds \right)^{p/4} \biggr).
	\end{align}
	By Holder's inequality, Young's inequality and A\ref{assump: mon, lip} one can bound the fifth term of \eqref{eq:expect_sup_estimate} as
	\begin{align} \label{eq:integral_estimate}
	E ^x \left( \int^t_0 \lvert \sigma(X_n(\kappa_n(s))) X_n(s)  \rvert^2 ds \right)^{p/4} & \leq C  \int^t_0 E ^x\lvert \sigma(X_n(\kappa_n(s))) X_n(s)  \rvert^{p/2} ds \nonumber \\
	& \leq C \left( 1+ \int^t_0 \left(E ^x\lvert X_n(\kappa_n(s)) \rvert^p+ E ^x \lvert X_n(s) \rvert^p \right) ds \right) .
	\end{align}
	Moreover, by Proposition \ref{th: strong mon} (ii), Young's inequality and
	\begin{align}
	\lvert b_n(X_n(\kappa_n(s))) &(X_n(s) - X_n(\kappa_n(s))) \rvert \nonumber \\
	& \leq n^{-1}\lvert b_n(X_n(\kappa_n(s))) \rvert^2+  b_n(X_n(\kappa_n(s))) \sigma(X_n(\kappa_n(s)) ) (W_s-W_{\kappa_n(s)}) \nonumber \\
	& \leq C \left( 1+ \lvert  X_n(\kappa_n(s)) \rvert^2 + n^{1/2} (1+\lvert  X_n(\kappa_n(s)) \rvert^2) \lvert W_s-W_{\kappa_n(s)}\rvert \right), \nonumber
	\end{align}
	one obtains, due to the independence of $X_n(\kappa_n(s))$ and $W_s-W_{\kappa_n(s)}$, that
	\begin{align} \label{eq:b_n_one-step_error_estimate}
	E^x \lvert b_n(X_n(\kappa_n(s))) (X_n(s) - X_n(\kappa_n(s))) \rvert^{p/2}\leq C \left( 1+ E^x \lvert  X_n(\kappa_n(s)) \rvert^p \right).
	\end{align}
	Then, substituting \eqref{eq:integral_estimate} and \eqref{eq:b_n_one-step_error_estimate} into \eqref{eq:expect_sup_estimate} yields
	\begin{align}
	E ^x \sup_{0 \leq t \leq T} \lvert X_n(t) \rvert^p &\leq C \biggr (1+ \lvert x\rvert ^p +\int^T_0 E ^x\lvert X_n(s)\rvert^p  + E ^x\lvert X_n(\kappa_n(s))\rvert^pds \biggr) \nonumber \\
	& \leq C \biggr (1+ \lvert x\rvert ^p +\int^T_0 E ^x\sup_{0\leq u \leq s}\lvert X_n(u)\rvert^p  ds \biggr)    , \nonumber
	\end{align}
	at which point applying Gronwall's inequality results in
	\begin{align}
	E ^x \sup_{0 \leq t \leq T} \lvert X_n(t) \rvert^p &\leq C(1+\lvert x\rvert ^p). \nonumber
	\end{align}
	Moreover, for any $0< q < 4$, raising to the power of $q/p$
	\begin{align}
	E ^x \sup_{0 \leq t \leq T} \lvert X_n(t) \rvert^q \leq (E ^x \sup_{0 \leq t \leq T} \lvert X_n(t) \rvert^p)^{q/p} &\leq C(1+\lvert x\rvert ^q). \nonumber
	\end{align}
	Furthermore since C does not depend on $x$, and since $\eta$ is independent of the driving noise, one observes that for every $0 \leq p \leq p_0$
	\begin{align}
	E \sup_{0 \leq t \leq T} \lvert X_n(t) \rvert^p \leq C(1+E\lvert \eta \rvert^p ).
	\end{align}
\end{proof}

\begin{Lemma} \label{lem: b term 2}
	Let $p(l+1) \leq p_0$. Then, for every $n \geq 1$
	\begin{equation}
	E \lvert b(X_n(t))- b(X_n(\kappa_n(t))) \rvert^p \leq Cn^{-p/2}. \nonumber
	\end{equation}
\end{Lemma}
\begin{proof}
	By A\ref{assump: poly growth} and Young's inequality one obtains that
	\begin{align}\label{eq: one-step comp b}
	\lvert b(X_n(t))- &b(X_n(\kappa_n(t))) \rvert^p \leq C(1+ \lvert X_n(t) \rvert^{pl} + \lvert X_n(\kappa_n(t)) \rvert^{pl}) \lvert  X_n(t)- X_n(\kappa_n(t))  \rvert^p \nonumber \\
	& \leq C((1+ \lvert X_n(\kappa_n(t)) \rvert^{pl}) \lvert  X_n(t)- X_n(\kappa_n(t))  \rvert^p+ \lvert  X_n(t)- X_n(\kappa_n(t))  \rvert^{(l+1)p} ),
	\end{align}
	where the inequality  $\lvert X_n(t) \rvert^{pl} \leq C(\lvert X_n(\kappa_n(t)) \rvert^{pl} + \lvert  X_n(t)- X_n(\kappa_n(t))  \rvert^{pl} )$ is used. Then, for every $0<q<p_0$ and in view of Proposition \ref{th: strong mon} (ii),  one obtains
	\begin{align} \label{eq: one step dif}
	E [ \lvert  X_n(t)- X_n(\kappa_n(t))  \rvert^q \vert \mathcal{F}_{\kappa_n(t)} ]&\leq C( n^{-q}  \lvert b_n(X_n(\kappa_n(t)))\rvert^q + \lvert \sigma(X_n(\kappa_n(t)))\rvert^q E\lvert W_t-W_{\kappa_n(t)} \rvert^q) \nonumber\\
	& \leq Cn^{-q/2}( 1+ \lvert X_n(\kappa_n(t) \rvert^q ).
	\end{align}
	Thus, by taking the conditional expectation of \eqref{eq: one-step comp b} with respect to $\mathcal{F}_{\kappa_n(t)}$  as in \eqref{eq: one step dif} and by using the estimate from \eqref{eq: one step dif}, one obtains upon the application of expectations
	\begin{align}
	E \lvert b(X_n(t))- &b(X_n(\kappa_n(t))) \rvert^p \leq Cn^{-p/2}(1+ E\lvert X_n(\kappa_n(t) \rvert^{p(l+1)} ), \nonumber
	\end{align}
	and consequently the result follows from Proposition \ref{thm: moment bounds}.
\end{proof}

\begin{proof}
	\textbf{Proof of Theorem \ref{thm:rate_of_convergence}} Firstly let us define
	\begin{equation}
	e_n(t) = X_n(t)-X(t). \nonumber
	\end{equation}
	Let us first assume $p \geq 4$. Then
	\begin{align}
	E \sup_{0\leq t \leq T} \lvert e_n(t) \rvert^p\leq CE \sup_{0\leq t \leq T} ( \lvert X_n(t) \rvert^p+ \lvert X(t) \rvert^p) < \infty,
	\end{align}
	and by It{\^o}'s lemma
	\begin{align}\label{eq: e_n^2}
	\lvert &e_n(t) \rvert^2 \leq C \biggr ( \int^t_0 (b(X(s)) -b_n(X_n(\kappa_n(s))) ) e_n(s) ds \nonumber \\
	& +\int^t_0 \lvert \sigma(X(s))- \sigma(X_n(\kappa_n(s))) \rvert^2 ds + \int^t_0 ( \sigma(X(s))- \sigma(X_n(\kappa_n(s))) ) e_n(s) dW_s \biggr ).
	\end{align}
	For the first term on the right hand side of \eqref{eq: e_n^2} one uses the splitting
	\begin{align}\label{eq: splitting}
	(b(X(t))-b_n(X_n(\kappa_n(t)))& e_n(t) = B_{1,n}(t)+B_{2,n}(t)+B_{3,n}(t),
	\end{align}
	for
	\begin{align}
	B_{1,n}(t):= ( b(X(t)) - b(X_n(t)))e_n(t) \leq -L \lvert e_n(t) \rvert^2, \nonumber
	\end{align}
	\begin{align}
	B_{2,n}(t) &:= ( b(X_n(t))- b(X_n(\kappa_n(t))) )e_n(t) \nonumber \\
	& \leq C( \lvert e_n(t) \rvert^2 +  \lvert b(X_n(t))- b(X_n(\kappa_n(t))) \rvert^2), \nonumber
	\end{align}
	\begin{align}
	B_{3,n}(t) &:= ( b(X_n(\kappa_n(t))) -b_n(X_n(\kappa_n(t))) )e_n(t) \nonumber \\
	& \leq C(  \lvert b(X_n(\kappa_n(t))) -b_n(X_n(\kappa_n(t))) \rvert^2 +\lvert e_n(t) \rvert^2 ), \nonumber
	\end{align}
	by A\ref{assump: mon, lip} and Young's inequality. Then
	\begin{align}
	(b(X(t))-b_n(X_n(\kappa_n(t)))& e_n(t) \leq C(\lvert e_n(t) \rvert^2 +  \lvert b(X_n(t))- b(X_n(\kappa_n(t))) \rvert^2 \nonumber \\
	& + \lvert b(X_n(\kappa_n(t))) -b_n(X_n(\kappa_n(t))) \rvert^2  ),
	\end{align}
	and substituting this into \eqref{eq: e_n^2}, raising to the power $p/2$ and applying Holder's inequality one obtains
	\begin{align}
	\lvert e_n(t) \rvert^p &\leq C \biggr ( \int^t_0\lvert e_n(s) \rvert^p +  \lvert b(X_n(s))- b(X_n(\kappa_n(s))) \rvert^p + \lvert b(X_n(\kappa_n(s))) -b_n(X_n(\kappa_n(s))) \rvert^p ds \nonumber \\
	& +\int^t_0 \lvert \sigma(X(s))- \sigma(X_n(\kappa_n(s))) \rvert^p ds + \biggr (\int^t_0 ( \sigma(X(s))- \sigma(X_n(\kappa_n(s))) ) e_n(s) dW_s\biggr ) ^{p/2} \biggr ), \nonumber
	\end{align}
	and therefore by Davis-Burkholder-Gundy's inequality and Holder's inequality
	\begin{align}\label{eq: exp of sup}
	E\sup_{0\leq t \leq T} \lvert e_n(s)& \rvert^p \leq C \biggr ( \int^T_0 E\lvert e_n(s) \rvert^p + E \lvert b(X_n(s))- b(X_n(\kappa_n(s))) \rvert^p  ds \nonumber \\
	&+\int^T_0 E\lvert b(X_n(\kappa_n(s))) -b_n(X_n(\kappa_n(s))) \rvert^p ds \nonumber \\
	& +\int^t_0 [E \lvert \sigma(X(s))- \sigma(X_n(\kappa_n(s))) \rvert^p +  E\lvert( \sigma(X(s))- \sigma(X_n(\kappa_n(s))))  e_n(s)\rvert^{p/2}] ds \biggr ).
	\end{align}
	Furthermore by Young's inequality again
	\begin{align} \label{eq: sigma bound 1}
	\lvert( \sigma(X(t))- \sigma(X_n(\kappa_n(t)))  e_n(t)\rvert^{p/2} \leq C( \lvert \sigma(X(t))- \sigma(X_n(\kappa_n(t))) \rvert^p+\lvert e_n(t)\rvert^p ),
	\end{align}
	and by A\ref{assump: mon, lip} and \eqref{eq: one step dif}
	\begin{align}\label{eq: sigma bound 2}
	E\lvert \sigma(X(t))- \sigma(X_n(\kappa_n(t))) \rvert^p &\leq CE\lvert X(t)- X_n(\kappa_n(t)) \rvert^p \nonumber \\
	&\leq C(E\lvert e_n(t) \rvert^p+  E\lvert X_n(t)- X_n(\kappa_n(t)) \rvert^p) \nonumber \\
	& \leq C(E\lvert e_n(t) \rvert^p+n^{-p/2}E\lvert X_n(t) \rvert^p ).
	\end{align}
	Then one can substitute \eqref{eq: sigma bound 1} and \eqref{eq: sigma bound 2} into \eqref{eq: exp of sup}, at which point applying Proposition \ref{th: strong mon} and Lemma \ref{lem: b term 2}
	\begin{align}
	E\sup_{0\leq t \leq T} \lvert e_n(t) \rvert^p &\leq C \biggr (\int^T_0 E\lvert b(X_n(t))- b(X_n(\kappa_n(t))) \rvert^p+ E\lvert b(X_n(\kappa_n(t))) -b_n(X_n(\kappa_n(t))) \rvert^p dt \nonumber \\
	&+ n^{-p/2}E\lvert X_n(t) \rvert^p+ \int^T_0 E\sup_{0\leq u \leq s} \lvert e_n(u) \rvert^p ds  \biggr ) \nonumber\\
	&\leq C \biggr (n^{-p/2}+n^{\frac{1}{2(l+1)} ( (l+2)p-p_0)}+ \int^T_0 E\sup_{0\leq u \leq s} \lvert e_n(u) \rvert^p ds  \biggr ) , \nonumber
	\end{align}
	so by Gronwall's inequality
	\begin{align}\label{eq: gronwell for dif}
	E\sup_{0\leq t \leq T} \lvert e_n(t) \rvert^p \leq C (n^{-p/2}+n^{\frac{1}{2(l+1)} ( (l+2)p-p_0)}).
	\end{align}
	For $0 \leq q < 4$ we can set $p=4$ in \eqref{eq: gronwell for dif} and raise both sides to the power of $q/4$. Then since	$E\sup_{0\leq t \leq T} \lvert e_n(t) \rvert^q \leq( E\sup_{0\leq t \leq T} \lvert e_n(t) \rvert^4 )^{q/4}$ the result follows.
\end{proof}

\section{Improved Rate for Constant Diffusion}
	
In the constant diffusion case (assuming some smoothness of $b$) one can achieve an improved rate of $n^{-p}$ in Theorem \ref{thm:rate_of_convergence}. Observing that in the constant diffusion case Milstein schemes and Euler schemes coincide, we show this using a similar argument to \cite{kumar2016milstein}. Firstly consider the following

\begin{assumptionA}  \label{assump: const diff} The function $b$ is of class $C^1$, and there exist positive constants $S>0$, $\sigma_0 \in \mathbb{R}^{d \times m}$, $l\geq1$ for which
\begin{equation}
\lvert Db(x)-Db(y) \rvert \leq S(1+\lvert x\rvert +\lvert y \rvert)^{l-1}\lvert x-y \rvert, \;\;\; x,y \in \mathbb{R}^d,
\end{equation}
\begin{equation}
\sigma(x) :=\sigma_0,
\end{equation}
where $Db$ is the Jacobian matrix of $b$.
\end{assumptionA}

\begin{Note}
	A\ref{assump: const diff} is a strengthening of A\ref{assump: poly growth} up to a change of constant, so we can apply the Lemmas from the previous sections by assuming A\ref{assump: const diff} instead of A\ref{assump: poly growth}.
\end{Note}

\begin{Th}\label{thm: conv const diff}
	Let A\ref{assump: mon, lip}, A\ref{assump: in cond} and A\ref{assump: const diff} hold. Then for every $p < p_0/(l+2)$ and $n \geq m$
	\begin{equation}
	\sup_{0\leq t \leq T}E[\lvert X(t) - X_n(t) \rvert^p] \leq Cn^{-r}, \nonumber
	\end{equation}
	where $r= \min \{\frac{1}{2(l+1)}(p_0(1\vee p/4)-p), p \}$.
\end{Th}

\begin{Lemma} \label{lem: Db bound}
	Let A\ref{assump: const diff} hold and suppose $p_0 \geq p(l-1)$. Then for every $x,y \in \mathbb{R}^d$
	\begin{equation}
	E\lvert b(X_n(s)) -b(X_n(\kappa_n(s))) -Db(x)(X_n(s)-X_n(\kappa_n(s)))\rvert^p \leq Cn^{-p}. \nonumber
	\end{equation}
\end{Lemma}
\begin{proof}
	Using \cite{kumar2016milstein}, Lemma 5, one estimates
	\begin{align}
	\lvert b(x) -b(y) -Db(x)(x-y) \rvert &\leq C(1+\lvert y \rvert +\lvert x \rvert)^{l-1}\rvert x-y \rvert^2 \\ \nonumber
	&\leq C(1+\lvert y \rvert^{l-1}+\lvert x-y \rvert^{l-1})\lvert x-y \rvert^{2}, \nonumber
	\end{align}
	so that one can consider $(1+\lvert y \rvert^{l-1})\lvert x-y \rvert^{2}$ and $\lvert x-y \rvert^{p+1}$ separately to obtain via Proposition \ref{thm: moment bounds} and Proposition \ref{th: strong mon} ii) that
	\begin{align}
	E \lvert b(X_n(s)) &-b(X_n(\kappa_n(s))) -Db(x)(X_n(s)-X_n(\kappa_n(s)))\rvert^p  \nonumber \\
	& \leq C(n^{-2p} (1+E\lvert X_n(\kappa_n(s)) \rvert^{p(l-1)} )\lvert b_n(X_n(\kappa_n(s)))\rvert^{2p} \nonumber \\
	&+(1+E\lvert X_n(\kappa_n(s)) \rvert^{p(l-1)}) \lvert \sigma_0 \rvert^{2p}E(W_s-W_{\kappa_n(s)})^{2p}  +n^{-p(l+1)/2}) \nonumber \\
	&\leq Cn^{-p}.
	\end{align}

\end{proof}
\begin{proof}
	\textbf{Theorem \ref{thm: conv const diff}} Let $p\geq 2$. Let $e_n(t)$ be as in the proof of Theorem \ref{thm:rate_of_convergence}. Then since $e_n(t)$ has vanishing diffusion coefficient, one can use the splitting from \eqref{eq: splitting} to obtain via A\ref{assump: mon, lip}, Proposition \ref{th: strong mon} and Lemma \ref{lem: b term 2} that
	\begin{align} \label{eq: e_n bound constdiff}
	E	\lvert e_n(t) \rvert^p& \leq C\int^t_0 E\lvert e_n(s) \rvert^{p-2} e_n(s) b(X(s))- b_n(X_n(\kappa_n(s))) ds \nonumber \\
	& \leq C\biggr(\int^t_0 \sup_{0 \leq u \leq s}E\lvert e_n(u) \rvert^{p} ds+ n^{\frac{1}{2(l+1)} ( (l+2)p-p_0)} +n^{-p} \nonumber \\
	&+ \int^t_0  E\lvert e_n(s) \rvert^{p-2} e_n(s) [ b(X_n(s))- b(X_n(\kappa_n(s)))] ds \biggr ),
	\end{align}
	One then splits the last term as
	\begin{align}\label{eq: calc start}
	& E\lvert e_n(s) \rvert^{p-2} e_n(s) [ b(X_n(s))- b(X_n(\kappa_n(s)))] \nonumber \\
	&=(\lvert e_n(s) \rvert^{p-2} e_n(s) - \lvert e_n(\kappa_n(s)) \rvert^{p-2} e_n(\kappa_n(s)) )[ b(X_n(s))- b(X_n(\kappa_n(s)))] \nonumber \\
	&+ \lvert e_n(\kappa_n(s)) \rvert^{p-2} e_n(\kappa_n(s)) [ b(X_n(s))- b(X_n(\kappa_n(s)))] :=I_{1,n}(s)+I_{2,n}(s).
	\end{align}
	To control $I_{1,n}$ one first observes
	\begin{equation*}
	\lvert \lvert e_n(s)  \rvert^{p-2} e_n(s)  - \lvert e_n(\kappa_n(s)) \rvert^{p-2} e_n(\kappa_n(s)) \rvert \leq C \sum_{i=1}^d \lvert \lvert e_n(s)  \rvert^{p-2} e^i_n(s)  - \lvert e_n(\kappa_n(s)) \rvert^{p-2} e^i_n(\kappa_n(s)) \rvert,
	\end{equation*}
	 and It{\^o}'s formula applied to $x \mapsto \lvert x \rvert ^{p-2}x^i$. For the latter, one uses additionally that the function $v_i(x): =\lvert x \rvert^{p-2} x^i$ satisfies $\lvert \nabla v_i(x)\rvert \leq C \lvert x \rvert^{p-2}$ for $i=1,2,...,d$. This then yields an expression for $\lvert \lvert e_n(s) \rvert^{p-2} e_n(s) - \lvert e_n(\kappa_n(s)) \rvert^{p-2} e_n(\kappa_n(s)) \rvert$. Then applying the splitting \eqref{eq: splitting} again, as well as Young's inequality and A\ref{assump: poly growth} (which follows from A\ref{assump: const diff}), one obtains
	\begin{align}
	I_{1,n}&(s) \leq CE \lvert b(X_n(s))- b(X_n(\kappa_n(s))) \rvert \int^s_{\kappa_n(s)} \lvert e_n(u) \rvert^{p-2} \lvert b(X(u))- b_n(X_n(\kappa_n(u)) \rvert  du \nonumber \\
	&\leq CE \lvert b(X_n(s))- b(X_n(\kappa_n(s))) \rvert \biggr ( \int^s_{\kappa_n(s)} \lvert e_n(u) \rvert^{p-1}(1+\lvert X(u) \rvert+\lvert X_n(u) \rvert)^l  du \nonumber \\
	&+ \int^s_{\kappa_n(s)} (\lvert e_n(u) \rvert^{p-2} \lvert b(X_n(u))- b(X_n(\kappa_n(u))) \rvert + \lvert e_n(u) \rvert^{p-2} \lvert b(X_n(\kappa_n(u))) - b_n(X_n(\kappa_n(u))) \rvert )du \biggr)  \nonumber \\
	&:=I_{1,1,n}(s)+I_{1,2,n}(s)+I_{1,3,n}(s).
	\end{align}
	Consequently, due to Theorem  \ref{thm: moment bounds} and the assumption $p_0 \geq p(2l+1)$, one obtains via Young's inequality
	\begin{align}
	I_{1,1,n}(s) &\leq C \int^s_{\kappa_n(s)} nE\lvert e_n(u) \rvert^{p}+n^{-p+1} E \lvert b(X_n(s))- b(X_n(\kappa_n(s))) \rvert ^p(1+\lvert X(u) \rvert+\lvert X_n(u) \rvert)^{pl} du \nonumber \\
	& \leq C (\sup_{0 \leq u \leq s} E\lvert e_n(u) \rvert^p + n^{-p}),
	\end{align}
	and furthermore by Lemma \ref{lem: b term 2}
	\begin{align}
	I_{1,2,n}(s) &\leq C \biggr ( \int^s_{\kappa_n(s)} ( nE\lvert e_n(u) \rvert^{p}+n^{-p/2+1} E \lvert b(X_n(s))- b(X_n(\kappa_n(s))) \rvert^p) du\nonumber \\
	&+n^{-p/2+1}\int^s_{\kappa_n(s)} E \lvert b(X_n(u)- b(X_n(\kappa_n(u))) \rvert^p  du \biggr )\nonumber \\
	& \leq C (\sup_{0 \leq u \leq s}E \lvert e_n(u) \rvert^p + n^{-p}),
	\end{align}
	and similarly for $I_{1,3,n}$, using Proposition \ref{th: strong mon} iii)
	\begin{align}
	I_{1,3,n}(s) &\leq C \biggr ( \int^s_{\kappa_n(s)} ( nE\lvert e_n(u) \rvert^{p}+n^{-p/2+1} E \lvert b(X_n(s))- b(X_n(\kappa_n(s))) \rvert^p) du\nonumber \\
	&+n^{-p/2+1}\int^s_{\kappa_n(s)} E \lvert b(X_n(\kappa_n(u)))- b_n(X_n(\kappa_n(u))) \rvert^p  du \biggr )\nonumber \\
	& \leq C (\sup_{0 \leq u \leq s}E \lvert e_n(u) \rvert^p + n^{-p}+n^{-\frac{p_0-p}{2(l+1)}}).
	\end{align}
	For $I_{2,n}$ we split as
	\begin{align}
	I_{2,n}(s)=	 E\lvert e_n(\kappa_n(s)) \rvert&^{p-2} e_n(\kappa_n(s)) [ b(X_n(s))- b(X_n(\kappa_n(s)))-Db(X_n(\kappa_n(s)))(X_n(s) - X_n(\kappa_n(s)))]\nonumber \\
	&+E\lvert e_n(\kappa_n(s)) \rvert^{p-2} e_n(\kappa_n(s)) Db(X_n(\kappa_n(s)))(X_n(s) - X_n(\kappa_n(s))) \nonumber \\
	&:=  I_{2,1,n}(s)+ I_{2,2,n}(s),
	\end{align}
	so that by Young's inequality and Lemma \ref{lem: Db bound}
	\begin{align}
	I_{2,1,n} =C (\sup_{0 \leq u \leq s}E \lvert e_n(u) \rvert^p + n^{-p}),
	\end{align}
	and by Proposition \ref{th: strong mon} and A\ref{assump: const diff}, plus the independence of $\lvert e_n(\kappa_n(s)) \rvert^{p-2} e_n(\kappa_n(s))$  with $\sigma_0 (W_t-W_{\kappa_n(s)})$ for $t \geq \kappa_n(s)$, one may conclude
	\begin{align}\label{eq: calc end}
	I_{2,2,n}(s) &= n^{-p} E\lvert  Db(X_n(\kappa_n(s))) b_n(X_n(\kappa_n(s))) \rvert^p +E\lvert e_n(\kappa_n(s)) \rvert^{p} \nonumber \\
	& \leq C \biggr ( n^{-p}E  \lvert Db(X_n(\kappa_n(s)))  \rvert^p \lvert b(X_n(\kappa_n(s))) - b_n(X_n(\kappa_n(s)))  \rvert^p  \nonumber \\
	& + n^{-p}E  \lvert Db(X_n(\kappa_n(s))) \lvert b(X_n(\kappa_n(s))) \rvert^p+  E\lvert e_n(\kappa_n(s)) \rvert^{p}\nonumber \biggr ) \\
	& \leq C (\sup_{0 \leq u \leq s}E \lvert e_n(u) \rvert^p + n^{-p}+n^{-\frac{p_0-p}{2(l+1)}}).
	\end{align}
	Then substituting \eqref{eq: calc start} - \eqref{eq: calc end} into \eqref{eq: e_n bound constdiff} and applying Gronwall's inequality yields
	\begin{equation}
\sup_{0\leq t \leq T}	E[\lvert X(t) - X_n(t) \rvert^p] \leq C(n^{-p}+n^{-\frac{p_0-p}{2(l+1)}}), \nonumber
	\end{equation}
	so that one obtains $L^q$ convergence, $q<2$, by setting $p=2$ and raising to the power of $2/q$.
\end{proof}

\bibliographystyle{amsplain}

\bibliography{NewRef}
\end{document}